\documentclass{article}
\usepackage{graphicx} 
\usepackage{hyperref}
\usepackage{amssymb}
\usepackage{amsopn}
\usepackage{amsmath}
\usepackage{amsfonts}
\usepackage{stmaryrd}
\usepackage{MnSymbol}
\usepackage{empheq}
\usepackage{tikz}
\usetikzlibrary{arrows}
\usepackage{etoolbox}
\AtBeginEnvironment{align}{\setcounter{equation}{0}}

\title{Set Theory is interpretable in Class Ordering Theory}
\author{Zuhair Al-Johar}
\date{October 2023}

\begin{document}

\maketitle

\section{Abstract}

Here it is shown that standard set theory can be interpreted in a theory about order. The ordering here is about non-extensional flat classes, i.e. classes that are not elements of classes. So, stipulating a nearly well order over all those classes coupled together with indexing that order by elements of those classes, thereby having those elements serve as ordinals; this together with infinity and a replacement like axiom would be shown to interpret $\sf ZFC$. Moreover, it is shown that a suitable version of this order theory is bi-interpretable with Morse-Kelley set theory augmented with a well ordering on classes.

\section{Introduction}

This theory align itself along attempts to formalize set theory in order theory. In a series of articles Takeuti had formalized the theory of ordinal numbers in $G'LC$ and $FLC$ both of which are written in higher order logic, and also in the ordinary first order logic with equality [Tak53,54,57,65]. He proved that a model of $\sf ZFC$ can be constructed in those theories. Koepke and Koerwien defined the theory of sets of ordinals $\sf SO$ in bi-sorted first order logic with equality, and proved its bi-interpretability with $\sf ZFC$ [KK05]. Here, a simpler axiomatic system, we'd name as $\sf COT$, is presented about flat classes of ordinals equipped with a near well-ordering on them together with an indexing function on that ordering which recursively sends bounded classes to ordinal indices. It would be shown here, that a version of such a theory is bi-interpretable with a more powerful system than $\sf ZFC$, that of More-Kelley set theory equipped with a well ordering on classes. The proof of bi-interpretability is vastly simpler than the case of $\sf SO$ with $\sf ZFC$. On the other hand, there is a growing interest in Class Choice principles extending $\sf MK$ set theory [G14], some of which prove to be hard questions [GHK21]. The point is that those choice principles are provable by adding a well-order on classes to $\sf MK$.

\section{Class Ordering Theory ``COT"}

\textbf{Language:}  Bi-sorted FOL with identity. First sort in lower case ranging over ordinals. Second sort in upper case ranging over classes of ordinals. \bigskip
 
\textbf{Primitives:} $\{=, \varepsilon, \prec , \iota\}$\bigskip

\textbf{Sorting}: $=$ denoting equality holding only between same sort objects; $\varepsilon$ standing for class element-hood relation holding only from first to second sort objects; $\prec$ denoting prior relation holding only for second sort objects; $\iota$ standing for the ordinal index of a class, a partial unary function holding only from second to first sort objects. \bigskip

\textbf{Axioms:} Bi-sorted ID axioms +:

 \begin{align*}  &\textbf{1. Comprehension: }\exists X \forall x \, (x \ \varepsilon \ X \leftrightarrow \phi), \text { if } X \text { is not free in }\phi \\ & \textbf{ Define: }  X \equiv Y \iff \forall z (z \ \varepsilon \ X \leftrightarrow z \ \varepsilon \ Y) \\ &\textbf { Define: } X \equiv \{x \mid \phi\} \iff \forall x \, (x \ \varepsilon \ X \leftrightarrow \phi)\\ &\textbf {2. Well Ordering} \prec \text{is a co-well-ordering on classes} \\&  \textbf{     -Transitive:} X \prec Y \prec  Z  \to X  \prec Z \\ & \textbf{ -Co-Connected:}   X \not\equiv Y \leftrightarrow [X \prec Y \lor Y \prec X], \\ & \textbf{ -Well-Founded:}   \phi(X) \to \exists M: \phi(M) \land \forall Y \, (\phi(Y) \to \neg Y  \prec M);  \text{ for every formula } \phi \\ &\textbf{ Define: } x < y \iff  \{x\} \prec \{y\} \\& \textbf { Define:} \lim X = \min l: \forall x \ \varepsilon \ X (x < l) \\ & \textbf{3. Respective: } \lim X \leq y \ \varepsilon \ Y \to X \prec Y \\&\textbf{ Define: } \operatorname {bnd}(X) \iff \exists l: l= \lim(X)  \\& \textbf{4. Indexing: }  \forall \operatorname {bnd} X \, \exists \iota(X): \iota(X)= \lim \{\iota(Y) \mid Y \prec X\}  \\ &\textbf{5. Infinity: } \exists l: \exists k < l \land \forall x < l \, \exists y \, ( x < y <  l)     \\&\textbf{6. Boundedness: }  \operatorname {bnd}(X) \land  X \operatorname \sim Y \to \operatorname {bnd} (Y) 
 \end{align*}

Where "$\sim$" stands for equinumerousity defined in the usual way after existence of a class bijection. \bigskip

As a notation we use  $\{ x \mid \phi \}$ not to stand for a particular class up to equality, but to stand for classes up to co-extensionality. So, it can be treated as a variable ranging over all classes having the same membership. So, when an expression $\psi(\{x \mid \phi\})$ is written, it denotes $\forall X: [\forall x \, (x \ \varepsilon \ X \leftrightarrow \phi)] \to \psi(X)$ \bigskip

Lemma I: $<$ is a well-ordering on ordinals \smallskip

Proof: The proofs of transitivity and well-foundedness for $<$ are straightforward, what needs to be proved is connectedness of $<$, that is: $$x \neq y \iff (x < y \lor y < x)$$  For left to right direction we have: $x \neq y \iff \{x\} \not \equiv \{y\} $ [Definitions of $\equiv$ and Singletons], thus by Co-connectedness, we have: $\{x\} \prec \{y\} \lor \{y\} \prec \{x\}$, thus $x < y \lor y < x$ [Definition of $<$]

For the opposite direction: $x < y \lor y < x \to \{x\} \prec \{y\} \lor \{y\} \prec \{x\}$, thus by Co-connectedness, we have: $\{x\} \not \equiv \{y\}$, thus $x \neq y$ [Definitions of $\equiv$ and Singletons] \bigskip

Lemma II: Transfinite recursion over ordinals

$\phi(0) \land \forall x (\forall y ( y < x \to \phi(y)) \to \phi(x)) \to \forall x: \phi(x) $ \smallskip

Proof: as usual, suppose there is an $x$ such that $\neg \phi(x)$, by well-foundedness of $<$, take the first such $x$, then it cannot be $0$, and so we must have $\phi(x)$ (antecedent of the above rule), a contradiction. \bigskip

Proposition I: every indexed class is bounded

$\iota(X) = x \to \operatorname {bnd}(X)$ \smallskip

Proof: Suppose there is an indexed unbounded class, let $X$ to be one of the $\prec$-minimal indexed unbounded classes, then every subclass $Y$ of a class $\{x \mid x < \iota(X)\}$ would be bounded and we'd have $\iota(Y) < \iota(X)$ [Respective and Indexing], take any class of all indices of those subclasses, that are not elements of the classes they index, and any such class would be bounded by $\iota(X)$ yet it won't have an index, contradicting axiom of Indexing. \bigskip

Proposition II: There is no last ordinal

$\not \exists l: \forall x \, (x \leq l)$ \smallskip

Proof: if there is a last ordinal $l$, then there is a last bounded class that $l$ indexes, i.e. $i^{-1}(l)$ is the last bounded class. Call a class \emph{elite} if its limit is $l$. Now, let $M$ be an $\prec$-minimal elite class; take any class $ K\equiv\{x \mid x < \iota(M)\}$, it is clear that any subclass of $K$ would be indexed by an element of $K$, since indexing ordinals $\geq \iota(M)$ are all reserved for indexing elite classes. Take a class of all elements of $K$ that are not elements of the classes they index, this would be a diagonal bounded class and it won't have an index, contradicting axiom of Indexing.\bigskip

Proposition III: every ordinal is an index

$\forall x \exists X: x= \iota (X)$ \smallskip

Proof: $0= \iota(\varnothing) $, let all ordinals $<l$ be indices, take any class $X$ of all ordinals $< l$ . Now, it cannot be the case that all subclasses of $X$ are indexed by ordinals prior to $l$, since if it is, then we can simply run the diagonal argument over it. So, there are subclasses of $X$ that are indexed with ordinals $ \geq l$, take an $\prec$ least of those subclasses, and that would have its index being $l$ [Indexing axiom]. Thus all ordinals are indices [Lemma II]. \bigskip

Lemma III: every index is higher than any element in the class it indexes

$y \ \varepsilon \ X \to \iota(X) > y$ \smallskip

Proof: $0$ is greater than all elements of any null class. Now if all ordinals $< l$ are greater than all elements in the classes they index, then $l$ must be so! Remember from proposition III, $l$ is the index of a $\prec$-least subclass of a class $\{x \mid x < l\}$ that has an index $\geq l$, and so all elements of that subclass are $<l$. Thus establishing the result. [Lemma II]

\section {COT interprets ZFC}

Now, define $\in$ as: $$ y \in x \iff   y \ \varepsilon \   \iota^{-1}(x) \iff \exists X: y \ \varepsilon \ X \land x= \iota ( X) $$. It'll be shown that $\sf ZFC$ would be interpreted over first sort world using this membership relation.\bigskip

Proof:

$\textbf{Extensionality: } \forall z( z \in x \leftrightarrow z \in y) \to x=y$

Proof: from the antecedent, we have $\iota^{-1}(x) \equiv \iota^{-1}(y)$. It can be proved that for any two co-extensional classes $A,B$ if we have a class $X \prec A$, then we must have $X \prec B$, the proof is easy, since if otherwise we must have $B \prec X$, and therefore $B \prec A$ [Transitivity of $\prec$], violating Co-connectedness of $\prec$. That said, then $\iota^{-1}(x) $ and $\iota{-1}(y)$ would have the same preceding classes, thus would be assigned the same limit by $\iota$, that is $x=y$. Since all ordinals are indices of classes [Proposition III], then there are no Ur-elements \bigskip

$\textbf{Regularity: } \exists x (x \in a) \to \exists y (y \in a \land \not \exists c: c \in y \land c \in a)$

Proof: let $\phi(z) \iff z \ \varepsilon \ \iota^{-1} (a)$, then per well-foundedness of $<$ there must exist an ordinal $m$ such that $\phi(m)$  and every ordinal $< m$ doesn't satisfy $\phi$, thus no $\varepsilon$-element of $\iota^{-1}(m)$ can satisfy $\phi$ [Lemma III], since the elements of $\iota^{-1} (a)$ are the ordinals that satisfy $\phi$, then we have $m \in a$ and $m$ is $\in$-disjoint of $a$ [Definition of $\in$] \bigskip

$\textbf{Pairing: } \forall a \forall b \exists x \forall y (y \in x \leftrightarrow y=a \lor y=b)$

Proof: every ordinals $a,b$, have a limit on the maximal one of them [Proposition II], thus every class $\{a,b\}$ is bounded, thus has the index $ \iota (\{a,b\}) $ which witness $x$ of the above formula. \bigskip

$\textbf{Union: } \forall a   \exists x \forall y (y \in x \leftrightarrow \exists z \in a (y \in z))$

Proof: from [Lemma III] and the definition of $\in$, every $a$ has all its $\in$-elements being prior (i.e. $<$) to it, that is $\forall y (y \in a  \to y < a)$, so all $\in$-members of $\in$-members of $a$ would also be $< a$, and any class of all of those would be bounded by $a$ itself, and so must have an ordinal index which would witness $x$ of the above formula. \bigskip

$\textbf{Power: } \forall a   \exists x \forall y (y \in x \leftrightarrow y \subseteq a)$

Proof: any class indexed by $a$ is bounded by $a$ itself [Lemma III], so all subclasses of an $\iota^{-1}(a) $ class are $\prec$ to every singleton class $\{a\}$ [Respective], so the indices of all those subclasses are $< \iota (\{a\})$ [Indexing], so any class of all those indices is bounded, thus has an ordinal index and that index would witness $x$ of the above formula. \bigskip

$\textbf{Separation: } \forall a   \exists x \forall y (y \in x \leftrightarrow y \in a \land \phi)$; if $x$ doesn't occur free in $\phi$. \smallskip

Proof: Since sets here are ordinal indices of bounded classes, clearly any sublcass of a bounded class is bounded, and so has an index that would witness $x$ in the above formula. \bigskip

$\textbf{Infinity: }$ The index of a class of all ordinals prior to the first ordinal limit is an infinite set. \bigskip

$\textbf{Choice: }$ This follows trivially actually the relation $<$ is a global well-order over all ordinals, and so over all sets \bigskip

Up till now, all axioms of $\sf ZC$ are interpreted. \bigskip

$\textbf{Replacement: }$ This follows trivially from the axiom of Boundedness and the definition of $\in$

\section {Interpreting COT in MK}

By MK we mean Morse-Kelley set theory, formalized as $\sf NBG$ but with allowing class variables freely in comprehension [Men97]. However, we'll adopt what is basically the original Kelley's system [Kel75] which is the well known equivalent approach [Git14] of replacing the limitation of size axiom by two axioms, one asserting that the range of every function from a set domain is a set, the other is Global Choice. \bigskip

Working in ``$\sf ZFC \ + $$ \ V=L \ + \text{ there is an inaccessible }$" take $L_\kappa$ where $\kappa$ is the first inaccessible. Now $ \mathcal P (L_\kappa)$ is a model of $\sf MK$ and moreover there is a well ordering $R$ over $\mathcal P(L_\kappa)$. Elements of $L_\kappa$ are the $\operatorname {Sets}$ of $\sf MK$, while the rest of $\mathcal P(L_\kappa)$ are the proper classes of $\sf MK$. Now we define a \textbf{\textit{rank respective}} well ordering $\prec$ over $\mathcal P(L_\kappa)$, as: $$X \prec Y \iff \rho X < \rho Y \land (\rho X = \rho Y \to X \ R \ Y)$$; where $\rho$ is the usual rank function.  \bigskip

Define: $ {\sf Ord } = \{x \in L_\kappa \mid x \text { is von Neumann ordinal} \}$

Let $\Omega = \{\{0\}\}$

Define: $X$ is  $class^*$ if and only if $\Omega \in X \land \forall Y \in X (Y \neq \Omega \to Y \in {\sf ORD})$.

Define $\in^*$ as: $X \in^* Y \iff class^*(Y) \land X \in Y \land X \neq \Omega $.

Define $\prec^*$ as: $ X \prec^* Y \iff \{Z: Z \in^* X\} \prec \{Z: Z \in^* Y\}$

Let, ${\sf Classes^*} = \{ K \mid class^*(K)\}$ \bigskip

We'll show that $({\sf Classes^*}, {\sf Ord} , =, \in^*, \prec^*, \iota^*)$ is a model of $\sf COT$ \bigskip

\textbf{Comprehension}: $\exists X \in {\sf Classes^*} \forall x \in {\sf Ord}: x \in^* X \leftrightarrow \phi(x)$ 

Proof: for any $\phi$ we have the class $A= \{x \in {\sf Ord} \mid \phi(x)\}$, 

then easily get $X= A \bigcup \{\Omega \}$ \bigskip

\textbf{Well-Ordering:} $\prec^*$ is a well ordering on ${\sf Classes^*}$

Proof: Clearly the formula $X= \{ Y \mid Y \in^* X^*\} $  describes a one-to-one correspondence between $\mathcal P({\sf Ord})$ and $\sf Classes^*$, and since $\prec$ is a well-ordering on the former, then from the definition of $\prec^*$ it must be a well-ordering the later. \bigskip

Define:  $x <^* y \iff  \{x\} \prec^* \{y\}$

Define: $ \lim^* X = \min l: \forall x \in^* X (x <^* l)  $ \bigskip

\textbf{Respective:} $X,Y \in \mathcal {\sf Classes^*} \land \lim^* X \leq y \in^* Y \to X \prec^* Y $

Proof: No distinct ordinals have the same rank! Therefore, for any ordinals $a,b$, the sets $\{a\}; \{b\}$ are of distinct ranks. Moreover the rank strict smaller than relation "$<^\rho$" is a strict well-order on $\sf Ord$ and thus also on $\mathcal P_1({\sf Ord})$. That said, then for every $\sf Class^*$ $X$ we have $\lim^* X$ having strictly greater rank than all $\in^*$-elements of $X$, and so $\rho \{x \mid x \in^* X \} < \rho \{\lim^* X\}  $, accordingly if there is $\lim^* X \leq y \in^* Y$, then $\rho  \{x \mid x \in^* X \} < \rho \{y \mid y \in^* Y\} $ and so $\{x \mid x \in^* X \} \prec \{y \mid y \in^* Y\}$ by definition, and so $X \prec^* Y$ by definition\bigskip

Define: $\operatorname {bnd}(X) \iff X \in \mathcal P({\sf Ord}) \land \exists l: l=\lim X$ \bigskip

Theorem 1: $\operatorname {bnd}(X) \iff X \in L_\kappa \cap \mathcal P ({\sf Ord}) $

Proof: if $X$ is a set of ordinals in $L_\kappa$ that is bounded upwardly by an ordinal in $L_\kappa$, then it has a rank in $L_\kappa$, and so it is an element of $L_\kappa$. For the opposite direction if $X \in L_\kappa \cap \mathcal P ({\sf Ord}$), then it has a rank in $L_\kappa$ and that rank is a limit on all elements of $X$ and so $\operatorname {bnd}(X)$ \bigskip

Define: $\operatorname {bnd^*}(X) \iff X \in {\sf Classes^*} \land \exists l: l=\lim^* X$ \bigskip

Theorem 2: $| \{X : \operatorname {bnd^*}(X) \} | = | L_\kappa \cap \mathcal P ({\sf Ord}) | $

Proof: Clearly the function $\eta(X)= \{Y \mid Y \in^* X\}$ is a bijection from bounded$^*$ classes   to  bounded classes and the latter are the elements of $L_\kappa \cap \mathcal P ({\sf Ord})$ [Theorem 1] \bigskip

\textbf{Indexing:} Define $\iota^*$ recursively on the bounded$^*$ classes  using the well ordering $\prec^*$ as:  $$\iota^*(X)=\lim \{ \iota^*(Y) \mid Y \prec^* X \}$$ \ \ \ \textbf{ Infinity:} There exists a limit ordinal. 

Proof:  $\omega$ is a limit ordinal with respect to $<$ \bigskip

\textbf{Boundedness:} $\operatorname {bnd^*}(X) \land Y \in {\sf Classes^*} \land  |\eta (X)| = |\eta(Y)| \to \operatorname {bnd^*}(Y)$

Proof: Since we have $\operatorname {bnd^*}(X)$ then $\eta(X) \in L_\kappa$, so by Replacement we get $\eta(Y) \in L_\kappa$, and so $\operatorname {bnd^*}(Y)$ by definition. \bigskip

\section{Bi-interpretability of COT and Set Theory}

We set a claim that Extensional Tri-sorted $\sf COT$ is bi-interpretable with $\sf MK$ extended with a well ordering $\prec$ on all classes, the latter theory to be denoted by $\sf MKCWO$. \bigskip

 Extensional Tri-sorted $\sf COT$, hereafter symbolized simply as $\sf COT$, adds on top of the above mentioned axioms of $\sf COT$ the axiom of Extensionality over all classes, and also adds a third sort that ranges over the first two sorts, we can use $\mathfrak {A,B,C,X,Y,Z,...}$ to denote the third sort. So, we axiomatize: \smallskip
$\\\textbf{Sorting axioms: } \\\forall {\mathfrak X}: \exists x \, ({\mathfrak X}=x) \lor \exists X \, ({\mathfrak X}=X) \\ \forall x \exists {\mathfrak X}: x = {\mathfrak X} \\ \forall X \exists {\mathfrak X}: X={\mathfrak X} \\ \forall x \forall X: x \neq X $  

\subsection{Definitions and Terminology:}

The main definitions and terminology follows that of Visser in \textit{Categories of Theories and Interpretations} [Vis06]. However, for the purpose of self inclusiveness they'd be concisely reintroduced here.   

A translation $\tau$ of signature $\Sigma$ in signature $\Theta$, denoted $\tau: \Sigma \longrightarrow \Theta$, is a pair $\langle \delta, F \rangle$, where $\delta$ is a formula in signature $\Theta$ with exactly one free variable, and $F$ is a function such that for each $n$-ary predicate (an $n$-ary relation, or $n-1$-ary function) $P$ in $\Sigma$ we have $F(P)$ [also written as $F_\tau(P)$] being an $n$-ary predicate defined by a formula in signature $\Theta$ having exactly $n$ free variables. 

Any formula $a_i$ written in $\Sigma$, its translation using $\tau$ is denoted by $a_i^\tau$ and it is the sentence resulting from merely replacing each element $P$ of the signature of $\Sigma$ occurring in $a_i$ by $F_\tau(P)$ and relativising all quantifiers in $a_i$ by $\delta_\tau$. For connivance each translated predicate $F_\tau(P)$ is written as $P_\tau$, also quantification is denoted by $\forall x: \delta_\tau; \exists x: \delta_\tau$, it is understood that $x$ is the sole free variable of $\delta_\tau$.

Interpreting a multi-sorted theory follows Visser and Friedman  [FV14], so when using translation $\tau$, this is done by assigning a domain formula $\delta^i_\tau$ for each sort $i$, then we require that for any $n$-ary predicate  $P$ in the signature of that theory that has $t_1,..,t_n$ type pattern, then the formula $A$ defining $P_\tau$ must assure that $\forall x_1,..,\forall x_n (P_\tau(x_1,..,x_n) \to \delta_\tau^{t_1}(x_1) \land ...\land \delta_\tau^{t_n}(x_n))$  \smallskip

An identity translation on signature $\Theta$, denoted by $\sf id_\Theta$, is $\sf id: \Theta \longrightarrow \Theta$ where $\delta_{\sf id} := \Huge{\top} $, and for every element $P$ of $\Theta$, $P_{\sf id} := P$

Any two translations $\tau, \nu$, the composition $\tau \circ \nu$, abbreviated as $\nu \tau$, is a translation with $\delta_{\nu \tau} := \delta_\tau \land (\delta_\nu)^\tau$ and $P_{\nu \tau } := (P_\nu)^\tau$

 An interpretation $I_\tau$ of theory $U$ written in $\Sigma$ in theory $V$ written in $\Theta$, denoted $I_\tau: U \longrightarrow V$, is a triplet $\langle U, \tau, V \rangle$ such that $\tau: \Sigma \longrightarrow \Theta$ and for every axiom $a_i$ of $U$ we have $V \vdash a_i^\tau$. 

 We use similar terminology to that of translations, so for an interpretation $K$ we have $\tau_K$ means the translation used in $K$, i.e. $K=\langle .,\tau,. \rangle$, now $\delta_K, E_K, P_K, A^K$ are actually $\delta_{\tau_K}; E_{\tau_K}, P_{\tau_K}, A^{\tau_K}$

 Similarly we can have identity interpretations and compositions of interpretations, so $I_{\sf id}: U \longrightarrow U$ is the triplet $\langle U, {\sf id_\Sigma}, U \rangle$, is the identity interpretation on theory $U$ in signature $\Sigma$.

 Also, for interpretations $I_\tau: U \longrightarrow W; I_\nu: W \longrightarrow V$, the composition interpretation $I_\nu \circ I_\tau: U \longrightarrow V$, is the triplet $ \langle U, \nu \circ \tau , V \rangle =\langle U, \tau \nu, V \rangle $

 $F$ is an i-isomorphism from interpretation K to interpretation M, where $K: U \longrightarrow V; M: U \longrightarrow V$, if and only if, $F$ is a formula in the language of $V$ having exactly two free variables, that satisfy all the following:

\begin{align} & V \vdash x F y \to (x: \delta_K \land y: \delta_M) \\ & V \vdash (x: \delta_K \land y: \delta_M \land x E_K x' F y' E_M y ) \to x F y \\ & V \vdash \forall x: \delta_K \exists y: \delta_M x F y \\ & V \vdash (x F y \land x F y') \to y E_M y' \\ & V \vdash \forall y: \delta_M \exists x: \delta_K x F y \\ & V \vdash (x F y \land x'Fy) \to x E_K x' \\ & V \vdash \vec{x} F \vec{y} \to (P_K\vec{x} \leftrightarrow P_M \vec{y})
\end{align}

Here $\vec{x} F \vec{y}$ abbreviates $x_0 F y_0 \land ...\land  x_{n-1} F y_{n-1}$, for appropriate $n$; and $E $ is the identity relation in the signature of $U$. \bigskip

To establish the bi-interpretability claim, we need to prove that there are interpretations  
\begin{align} &\sf I_\tau: COT \longrightarrow  MKCWO \\ & \sf I_\nu: MKCWO \longrightarrow COT \end{align}, 
such that:
\begin {align} & \sf I_\nu \circ I_\tau: COT \longrightarrow COT \text { is i-isomorphic to }\sf I_{id}: COT \longrightarrow COT \\ & \sf I_\tau \circ I_\nu: MKCWO \longrightarrow MKCWO \text { is i-isomorphic to } \sf I_{id}: MKCWO\longrightarrow MKCWO \end{align}

\subsection{Proof}

\subsubsection{ First direction }Let $\nu: \mathcal L(=,\in, \prec) \longrightarrow \mathcal L(=,\varepsilon,\prec, \iota)$ \bigskip

Define: $\delta_\nu := \exists x( \mathfrak X = x) \lor \neg \operatorname {bnd}( \mathfrak X) $ \bigskip

Define: $ \mathfrak X =_\nu \mathfrak X \iff  \mathfrak X = \mathfrak X $

Define: \!  \begin{align*} {\mathfrak A} \in_\nu {\mathfrak X} \iff \exists a: {\mathfrak A}= a \land   [&(\exists x: {\mathfrak X}=x \land a \ \varepsilon \ \iota^{-1}(x)) \lor \\ & (\exists X: {\mathfrak X} = X \land a \ \varepsilon \ X)] \end{align*}  

Define: 
\begin{align*} \mathfrak A \prec_\nu \mathfrak X \iff & \exists x \exists a: \mathfrak X=x \land \mathfrak A= a \land a < x \ \lor \\& \exists X \exists A: \mathfrak X=X \land \mathfrak A=A \land A \prec X \ \lor \\ & \exists X \exists a: \mathfrak X=X \land \mathfrak A=a  
\end{align*}

This way we get to interpret $\sf MKCWO$, this is straightforward for the set axioms of $\sf MK$, i.e. set Extensionality, Pairing, Union, Power, Separation, Replacement, and Infinity, all follow by interpretability of $\sf ZFC$. Now, foundation for classes clearly follows because all sets are well founded so classes of them would be so. What remains is Class Extensionality, Comprehension and Global Choice. For the first, it follows from  $\in_\nu$-co-Extensionality preserving the state of being bounded or not, thus preserving set-hood and proper class-hood; then Class Extensionality straightforwardly follow from Extensionality of classes in $\sf COT$. Comprehension also follows easily from  Comprehension in $\sf COT$. Global Choice clearly follows from the well ordering on all classes. To be noted is that quantification over classes is allowed both in Comprehension as well as in Separation and Replacement.    

\subsubsection{ Second direction }

Let $\tau: \mathcal L(=,\varepsilon,\prec, \iota) \longrightarrow \mathcal L(=,\in, \prec)$ \smallskip

Without loss of generality we take $\prec$ to be rank respective [section 4]. \smallskip

Fix a non ordinal set $\Omega$; \smallskip
 
Define: 
$\\\delta^1_\tau:= x \in {\sf ORD} \\ \delta^2_\tau: = \Omega \in x \land \forall y \in x (y \neq \Omega \to y \in {\sf ORD}) \\ \delta^3_\tau: = \delta^1_\tau \lor \delta^2_\tau \\ x =_\tau y \iff x =y \\x \ \varepsilon _\tau \ y \iff  \delta^1_\tau(x) \land \delta^2_\tau(y) \land x \in y\\ \zeta(x) = i \iff \operatorname {Set}(x) \land i= \lim \operatorname {ord} \{\zeta(y) \mid y \prec x \}\\ x \prec_* y \iff x \subseteq {\sf ORD} \land y \subseteq {\sf ORD} \land \{\zeta^{-1}(m) \mid m \in x\} \prec \{\zeta^{-1}(n) \mid n  \in  y\} \\ x \prec_\tau y \iff \delta^2_\tau (x) \land \delta_\tau^2(y) \land \eta(x) \prec_* \eta(y)  \\\iota_\tau(x)= y \iff \delta^2_\tau(x) \land y=\zeta (\{ \zeta^{-1}(z) \mid z \ \varepsilon_\tau \ x \})$ \bigskip

Interpreting $\sf COT$ follows exactly the same argument in section 4. Just replace the domains appropriately, and switch $P^*$ to $P_\tau$. What needs to be proved is: $$\iota_\tau(x)= \lim \operatorname {ord} \{\iota_\tau(y) \mid y \prec_\tau x \} $$ Proof: it is obvious that the well ordering relation $\prec$ on sets is copied downwardly to   $\prec_*$ on $\mathcal P({\sf ORD})$, which is in turn copied downwardly to $\prec_\tau$ on $\sf Classes^*$ (i.e. $\delta^2_\tau$ objects). Then the $\zeta$ function which is an index of $\prec$ on sets, would be copied downwardly to the $\iota_\tau$ function on $\sf Classes^*$, and so $\iota_\tau$ would like-wisely index $\prec_\tau$ relation on $\sf Classes^*$. That said, it is easy to prove that whatever ordinal $\alpha$ we have $\alpha < \iota_\tau(x) \iff \exists y: \delta^2_\tau \ y \prec_\tau x \land \alpha=\iota_\tau(y)$. Thus establishing the result.

\subsubsection{ Proving i-isomorphisms }

\textbf{First direction:}

That $F$ is i-isomorphic from $K$ to $M$ where $K=\sf \langle MKCWO, {\sf id}, MKCWO \rangle $ and $M=\sf \langle MKCWO, \tau \circ  \nu,  MKCWO \rangle $  , is proved as follows:\bigskip

Working in $\sf MKCWO$, define: 
\begin{align*}
F (X) = 
\left\{
    \begin {aligned}
         & \zeta(X) \quad & \text {if } \operatorname {Set}(X) \\
         &  \{\zeta(Y) \mid Y \in X \} \cup \{\Omega\}\quad & \text{if }  \operatorname  {\neg Set}(X)                
    \end{aligned}
\right.
\end{align*}

$\delta_K = \delta_{\sf id} := x \in \mathcal P (L_\kappa)$

$\delta_M = \delta_{\nu\tau} :=  \delta^3_\tau \land  (\delta_\nu)^\tau $

$\delta_M:= y \in {\sf ORD} \cup \{x \cup \{\Omega\} \mid x \subseteq {\sf ORD} \land \neg \operatorname {Set}(x)\}$

$x E_K x' \iff x = x'$

$ y E_M y' \iff y=y'$

$x \in_K x' \iff x \in x'$

$ y \in_M y' \iff (y \in_\nu y')^\tau $

$y \in_M y' \iff (y' \in {\sf ORD} \land y \in \iota_\tau^{-1} (y')) \lor (\delta^2_\tau(y') \land \neg \operatorname {Set}(y') \land y \in y')$  

$x \prec_K x' \iff x \prec x'$

$ y \prec_M y' \iff (y  \prec_\nu y')^\tau$
\begin{align*}  x \prec_M  y \iff & x,y \in {\sf ORD} \land x < y \ \lor \\& \delta^2_\tau(x) \land \delta^2_\tau(y) \land \neg \operatorname {Set}(x) \land \neg \operatorname {Set}(y) \land x \prec_\tau y \ \lor \\ & x \in {\sf ORD} \land \delta^2_\tau(y) \land \neg \operatorname {Set}(y) 
\end{align*}

It's an easy check to see that all 1-7 requirements of i-isomorphisms are met!\bigskip

\textbf{Second direction:}

That $G$ is i-isomorphic from $K=\sf \langle COT, {\sf id}, COT \rangle $ to $M=\sf \langle COT, \nu \circ  \tau,  COT \rangle $  , is proved as follows: \bigskip

Working in $\sf COT$, define: $\zeta_\nu  := (\lim \operatorname {ord} \{\zeta(y) \mid y \prec x \})^\nu$
\begin{align*}
\text {Define: } \ \ \ G (\mathfrak X) = 
\left\{
    \begin {aligned}
         & \zeta_\nu(\mathfrak X) \quad & \text { if } \exists x: \mathfrak X=x \\
         & \iota (\{\zeta_\nu(\mathfrak A) \mid \mathfrak A \ \varepsilon \ \mathfrak X\} \cup \{\Omega\})\quad & \text { if } \exists X: \mathfrak X=X \land \operatorname {bnd}(X)\\
         &  \{\zeta_\nu(\mathfrak A) \mid \mathfrak A \ \varepsilon \ \mathfrak X \} \cup \{\Omega\} \quad & \text{ if } \exists X: \mathfrak X= X \land \neg \operatorname {bnd}(X)               
    \end{aligned}
\right.
\end{align*}

 We'll adopt the notation $x^i$ to stand for the $i^{th}$ sort, so $\phi(x^1); \phi(x^2); \phi(x^3)$ actually means $ \phi(x), \phi(X), \phi(\mathfrak X)$ \smallskip

 $\delta^i_K  = \delta^i_{\sf id} :=   x^i = x^i$

$\delta^i_M = \delta^i_{\tau\nu} :=  \delta_\nu \land  (\delta^i_\tau)^\nu $

$x^i E_K y^i \iff x^i =y^i$

$z^i E_M u^i  \iff z^i =u^i$

$x^{1,3} \ \varepsilon_K \ y^{2,3} \iff x^{1,3} \ \varepsilon \ y^{2,3}$

    $ z^{1,3} \ \varepsilon_M \ u^{2,3} \iff (z^{1,3} \ \varepsilon_\tau \ u^{2,3})^\nu $

$x^{2,3}  \prec_K y^{2,3}  \iff x^{2,3}  \prec y^{2,3}$

$ z^{2,3} \prec_M u^{2,3} \iff (z^{2,3}  \prec_\tau u^{2,3})^\nu$

$\iota_K(x^{2,3}) = y^1 \iff \iota(x^{2,3}) = y^1 $

$\iota_M(z^{2,3})=u^1 \iff (\iota_\tau(z^{2,3})=u^1)^\nu$\bigskip

It's an easy check to see that all 1-7 requirements of i-isomorphism are met! \smallskip

More demonstrative, this means that given whatever model $\mathcal M$ of $\sf MKCWO$ we can use translation $\tau$ to build a model $\mathcal C_\tau^{\mathcal M}$ of $\sf COT$ within it, and within $\mathcal C_\tau^{\mathcal M}$, via translation $\nu$, we can build a model $M_\nu^{\mathcal C_\tau^{\mathcal M}}$ of $\sf MKCWO$, and such that we can define within $\mathcal M$ an isomorphism after formula $F$ from $ \mathcal M $ to  $M_\nu^{\mathcal C_\tau^{\mathcal M}}$. And also given whatever model $\mathcal C$ of $\sf COT$ we can use translation $\nu$ to build a model $\mathcal M_\nu^{\mathcal C}$ of $\sf MKCWO$ within it, and within $ \mathcal M_\nu^{\mathcal C}$, via translation $\tau$, we can build a model $C_\tau^{ \mathcal M_\nu^{\mathcal C}}$ of $\sf COT$, and such that we can define within $\mathcal C$ an isomorphism after formula $G$ from $ \mathcal C $ to  $C_\tau^{\mathcal \mathcal M_\nu^{\mathcal C}}$.

\section {Conclusion}

The powerful hierarchical class theory $\sf MKCWO$ is reducible in structure to a theory about indexed ordering of bounded classes of ordinals. It needs to be warned that this bi-interpretability is only with the extensional form of class ordering theory. The non-extensional form of $\sf COT$ is most likely not bi-interpretable with $\sf MKCWO$, but it might be bi-interpretable with a non-extensional version of it. The point is that the massive hierarchical structure of such a rich theory as $\sf MKCWO$ is in some sense redundant, everything can be done in a single tier of membership, i.e. with flat classes, and moreover the system is simple and the proofs of interpretation results in it are also simple. This motivates this system as a strong and simple tool for proving interpretability results. 

\section{References:}
\begin{itemize}
 \bibitem[Tak54]{Takeuti:Construction of the set theory from the theory of ordinal numbers.}
G. Takeuti.
\newblock ``Construction of the set theory from the theory of ordinal numbers".
\newblock J. Math. Soc. Japan, 6:196–220, 1954.
\bibitem [Tak53] {0} ------. ``On the generalized logic calculus." Jap. J. Math., 25:39-96, 1953; Errata to `On a generalized logic calculus', Jap. J. Math., 24:149-156, 1954.
\bibitem [Tak57] {1} ------. ``On the theory of ordinal numbers." J. Math. Soc. Japan, 9:93-113, 1957.
\bibitem [Tak65] {2} ------. ``A Formalization of the Theory of Ordinal Numbers." J. Symbolic Logic, 30(3):295–317, 1965.
\bibitem [Vis06] {3} A. Visser. ``Categories of theories and interpretations." In Logic in Tehran, Lecture Notes in Logic, pp. 284-341. Cambridge: Cambridge University Press, 2006.
\bibitem  [FV14] {4} H. M. Friedman and A. Visser. ``When Bi-Interpretability Implies Synonymy." Logic Group Preprint Series, 320:1–19, 2014.
\bibitem  [KK05] {5} P. Koepke and M. Koerwien. ``The theory of sets of ordinals." arXiv preprint math/0502265, 2005.
\bibitem  [Men97] {6} E. Mendelson. ``An Introduction to Mathematical Logic (4th ed.)." London: Chapman and Hall/CRC, 1997. Pp. 225–86. ISBN 978-0-412-80830-2.
\bibitem [GHK21] {7} V. Gitman, J. D. Hamkins, and A. Karagila. ``Kelley-Morse set theory does not prove the class Fodor principle." Fundamenta Mathematicae, 254(2), 133-154, 2021.
\bibitem [Git14] {8} V. Gitman.``Variants of Kelley-Morse set theory".  

https://victoriagitman.github.io/research/2014/02/24/variants-of-kelley-morse-set-theory.html
\bibitem[G14]{10} V. Gitman. ``Choice schemes for Kelley-Morse set theory".

https://victoriagitman.github.io/talks/2014/07/30/choice-schemes-for-kelley-morse-set-theory.html

\bibitem [Kel75] {9 } J. L. Kelley, General topology. New York: Springer-Verlag, 1975, p. xiv+298.
\end{itemize}

\end{document}